\documentclass[preprint,review,12pt]{elsarticle}

\usepackage[utf8]{inputenc}
\usepackage{amsmath}
\usepackage{amssymb}
\usepackage{graphics}
\usepackage{graphicx}
\usepackage{epstopdf}
\usepackage{xcolor}
\usepackage{caption}
\usepackage{subcaption}

\newtheorem{theorem}{Theorem}

\journal{}

\begin{document}

\begin{frontmatter}




\title{Accurate benchmark results of Blasius boundary layer problem using a leaping Taylor's series that converges for all real values.}


\author[inst1]{Anil Lal S. \corref{cor1} }

\affiliation[inst1]{organization={Department of Mechanical Engineering},
            addressline={Amrita Vishwa Vidyapeetham}, 
            city={Amritapuri Campus},
            postcode={690525}, 
            state={Kerala},
            country={India}}
\cortext[cor1]{Corresponding author: Anil Lal S.,  anillals@am.amrita.edu}

\author[inst2]{Milin Martin}

\affiliation[inst2]{organization={Department of Mechanical Engineering},
            addressline={College of Engineering, Trivandrum}, 
            city={Sreekariyam},
            postcode={695016}, 
            state={Kerala},
            country={India}}

\begin{abstract}
Blasius boundary layer solution is a Maclaurin series expansion of the function \(f(\eta)\), which has convergence problems when evaluating for higher values of \(\eta\) due to a singularity present at \(\eta\approx-5.69\).  In this paper we are introducing an accurate solution to \(f(\eta)\) using Taylor's series expansions with progressively shifted centers of expansion(Leaping centers). Each series is solved as an IVP with the three initial values computed from solution of the previous series, so the gap between the centers of two consecutive expansions is selected from within the convergence disc of the first series.   The last series is formed such that it is convergent for a reasonable high \(\eta\) value  needed for implementing the  boundary condition at infinity. The present methodology uses Newton-Raphson method to compute the  value of the unknown initial condition viz. \(f''(0)\) in an iterative manner. Benchmark accurate results of different parameters of flat plate boundary layer with no slip and slip boundary conditions have been reported in this paper.
\end{abstract}



\begin{keyword}
Blasius function \sep Leaping Taylor's series \sep Slip flow \sep Newton-Raphson method
\end{keyword}

\end{frontmatter}


\section{Introduction}
Blasius  series solution \cite{Blasius1908,Blasius1912} is well known as the accurate solution for flat plate boundary layer flow. It is not only a remarkable achievement in the history of fluid dynamics but also a groundbreaking  mathematical exercise in which the non-linear differential equation $f''' + \frac{1}{2}ff'' = 0$ with the known conditions $f(0) = 0, f'(0) = 0, f'(\infty) = 1$, called the Blasius boundary layer equation is solved.   It is not an initial value problem because \(f''(0)\) is not known and this poses certain difficulties in determining its solution.  The function \(f(\eta)\) represents the value of stream function and is popularly known as the Blasius function. Blasius determined the solution in the form of  a power series expansion  \cite{Boyd2008} as,
\begin{equation}
    f(\eta) = \frac{\kappa}{2}{\eta}^2 - \frac{{\kappa}^2}{240}{\eta}^5 + \frac{11}{161280}{\kappa}^3{\eta}^8 - \frac{5}{4257792}{\kappa}^4{\eta}^{11} + ...  \label{bl}
\end{equation}
where $f''(0) = \kappa$. The above series with center of expansion at \(\eta=0\) is found to have convergence problems at certain values of \(\eta\) and due to this the value of \(f''(0)\) is often found with a lesser accuracy. Understanding the mathematical nature of Blasius function including convergence properties, obtaining better accurate solutions and extension to cases of slip boundary conditions etc, have prompted the researchers to further investigate and solve the problem  analytically  \cite{Boyd1998,Catal2012} and numerically \cite{Cortell2005,Erturk2008}. 

Along with analytical methods, there have been efforts to use finite-order accurate numerical techniques such as Runge-Kutta method  and shooting method for getting the unkknown initial condition etc., due to the ease of computation. The shooting method \cite{white2006viscous} is used extensively, such that it has literally become a part of classroom exercises. The shooting method is a technique for solving a BVP by reducing it into  an IVP. The method consists of solving IVP with different initial conditions until the boundary condition of the BVP is satisfied. Another technique used to deduce the initial condition is called T\"{o}pfer transformation ($g = {\lambda}^{-1/3}f$) \cite{Topfer1912, Fazio2009}. T\"{o}pfer transformation leaves the initial conditions unchanged  along with the introduction of an additional initial condition, leading to the computation of the Blasius function and its derivatives through an inverse transformation. A number of previous investigations have reported the use of the Runge-Kutta scheme along with the shooting method \cite{Cortell2005,FANG20083088} to determine the unknown initial condition of the problem. One of the historical literature on numerical method adopted for solving the Blasius equation is by Howarth \cite{Howarth1937} in which  the value of $f''(0)$, accurate upto five decimal places is reported. A recent study by Majid et al.\cite{Majid2017} uses a predictor-corrector two-point block method to solve the Blasius equation numerically. The idea revolves around the development of a technique to bypass the step of reducing the BVP to an IVP, thus it makes use of a multiple shooting method that converges faster . 

Previous studies on solutions to Blasius function have also been attempted using analytical or semi-analytical methods  with an aim to get the solution with higher accuracy and with out convergence difficulties. Adomian Decomposition Method (ADM) is a semi-analytical method that uses ''Adomain Polynomials'' to represent the non-linear portion of the function. Abbaoui et al. \cite{Abbaoui1995} have given a proof for the convergence of the series solutions obtained from Adomaian Decomposition Method. Wang \cite{WANG20041} has shown the application of ADM to obtain an approximate analytical solution to Blasius function. Hashim \cite{Hashim2006} presents a modified technique that uses a combined ADM and Pade approach to report the improvement in the accuracy of the solution. Boyd \cite{Boyd1997} reported that series solution alone is not enough and  combining it with the Pade approximation is an effective approach. Pade approximation provides very close estimate of the function by writing it in the form of a rational function of a given order. It is reported that Pade Approximation to  Blasius equation with $\eta = 6$ provides convergence to the solution with a relative error of about $0.15\%$ in predicting the accurate value of $f''(0)$. For $\eta >8$, the convergence is found to be slow and there is a need for higher values of Pade order (N) for getting a better solution.  A recent study by Andrianov et al. \cite{Andrianov2020} introduced the Pade approximation for solving the boundary layer problem and for the problem of rotating fluid of Ekman layer type.

Differential Transformation Method (DTM) is a systematic framework that gives transformations which becomes the coefficients of terms of the series expansion of a function or functions  presented in the form of  linear, non-linear and simultaneous ordinary differential equations~ \cite{LIENTSAI1998101, zhou1986differential, ames1968nonlinear, Lal2014}. It may be noted that Blasius has derived his series accurately in the same manner but without resorting to DTM framework.  Yu et. al. \cite{LIENTSAI1998101}, used DTM by converting the Blasius problem into a set of three IVPs. A sub-domain approach is found to overcome the convergence problems of the series solution and values accurate up to \(6\) decimal places are reported. Arikoglu and Ibrahim \cite{Arikoglu2005a} have presented an inner-outer matching solution for the Blasius problem  in which \(f''(0)\) is evaluated by equating two series expansions with one series expanded about \(\eta=0\) and the other one about a higher value of \(\eta=m\) of about \(10\). The matching of the solutions is done at a value of \(\eta\) in between \(4\) and \(5\) to determine \(f''(0)\) and the reported value is found to have a difference starting from the fourth decimal place onward.  Lal and Neeraj  \cite{Lal2014} have used  DTM along with the scaled transformation for converting the boundary value of Blasius function  to initial value and reported results accurate up to 21 decimal points with a larger radius of convergence.

In the present work, we report a method in which a set of leaping series solutions that converge for all positive real values of \(\eta\) are developed, with out resorting to approximations. Application of Newton-Raphson method enables determination of any unknown initial value in the problem from a suitable boundary condition, and still obtain Taylor's series solution. Currently used methods such as T\"{o}pfer transformation, truncated boundary approximation, iterative transformation method \cite{Topfer1912,Fazio2009} are expected to pose difficulties in deriving the initial conditions from boundary conditions for the class problems governed by simultaneous differential equations and equations with more number of unknown initial values. We report benchmark solutions to boundary layers with no-slip and slip flow over the wall surface,  nature of convergence of series expansions, unknown initial values and parameters of asymptotic variation with a greater accuracy in this paper.

\section{Differential Transform}
 The Taylor series expansion of a function $f(\eta)$ about $\eta = \eta_{i}$ is,
\begin{equation}
  f(\eta)  = \sum_{k=0}^{\infty} \frac{\left.f^{k}(\eta)\right|_{\eta=\eta_i}}{k!} (\eta-\eta_i)^k= \sum_{k=0}^{\infty}F(k)(\eta-\eta_i)^k
\end{equation}
Where,  $F(k) = \frac{\left.f^{k}(\eta)\right|_{\eta=\eta_i}}{k!}$,~ \(\forall~k\in \mathbb{W}\) is called the differential transform of $f(\eta)$  and \(f^{k}(\eta)\) is the \(k^{th}\) derivative of \(f\) with respect to \(\eta\). We denote differential transform operation of \(f(\eta)\) as \(DT(f(\eta))\) and the result of the transform as \(F(k)\). Let \(U(k)\) and \(V(k)\) are the differential transforms of the functions \(u(\eta)\) and \(v(\eta)\). Some useful theorems of differential transforms with \(U(k)\) and \(V(k)\) \cite{Arikoglu2005, Arikoglu2008, Konuralp2011} that becomes handy in finding the solutions of IVPs are given here.
\begin{theorem}
If $f(\eta) = u(\eta) \pm v(\eta)$, then $F(k) = U(k) \pm V(k)$
\end{theorem}
\begin{theorem}
If $f(\eta) = \alpha u(\eta)$, then $F(k) = \alpha U(k)$
\end{theorem}
\begin{theorem}
If $f(\eta) = \frac{du(\eta)}{d\eta}$, then $F(k) = (k+1) U(k+1)$
\end{theorem}
\begin{theorem}
If $f(\eta) = \frac{d^{m}u(\eta)}{d\eta^{m}}$, then $F(k) = (k+1)(k+2)...(k+3) U(k+m)$
\end{theorem}
\begin{theorem}
If $f(\eta) = u(\eta)v(\eta)$, then $F(k) = \sum_{l=0}^{k} V(l)U(k-l)$
\end{theorem}
Making use of these theorems the terms in the Blasius equation namely $f'''(\eta), f''(\eta)$ and $f(\eta)$ can be differential transformed as \begin{equation}
    DT\left[f'''(\eta)\right]= (k+1)(k+2)(k+3)F(k)
\end{equation}
\begin{equation}
    DT\left[f(\eta)f''(\eta)\right] = \sum_{l=0}^{k}(l+1)(l+2)F(l)F(k-l)
\end{equation}
By knowing $F(0)$, $F(1)$, the value of $F(2)$ taken based on generated values of $f''(0)$, the  differential transform  for any value of $k$ is found using
\begin{equation}
    F(k+3)=-\frac{\sum_{l=0}^{k}(l+1)(l+2)F(l)F(k-l)}{2(k+1)(k+2)(k+3)}
\end{equation}
For \(\eta_i=0\), we have \(F(0)=0\), \(F(1)=0\) and \(F(2)=\kappa/2\), where \(\kappa=f''(0)\).

\section{Solution Methodology}
\subsection{Series expansion about leaping centers}
It is well known that Taylor's series given in equation (1) with center of expansion at $\eta = 0$ does not converge for $\eta >5.69$ \cite{Boyd2008}. But in computations on computers, convergence up to \(\eta=5.69 \) is rarely possible due to accuracy problems and restricted summation over a finite number of terms. This leads to approximating the boundary condition  \(f'(\infty)=1\) as \(f'(5.69)=1\), which results in lesser accurate results.  In the present work, we propose an  approach of combining more number of series expansions of the Blasius function by leaping  the center of expansion to different values of $\eta$ to circumvent the convergence problems associated with approximation of  \(f'(\infty)=1\). We have approximated  the infinity  by a larger value of \(\eta\) and applied the boundary condition as \(f'(\eta)=1\), where \(\eta\ge18\). Such an approximation is found to provide a more  accurate solution to the Blasius function.  
\begin{figure}[htp!]
    \centering
    \includegraphics[scale = 0.55]{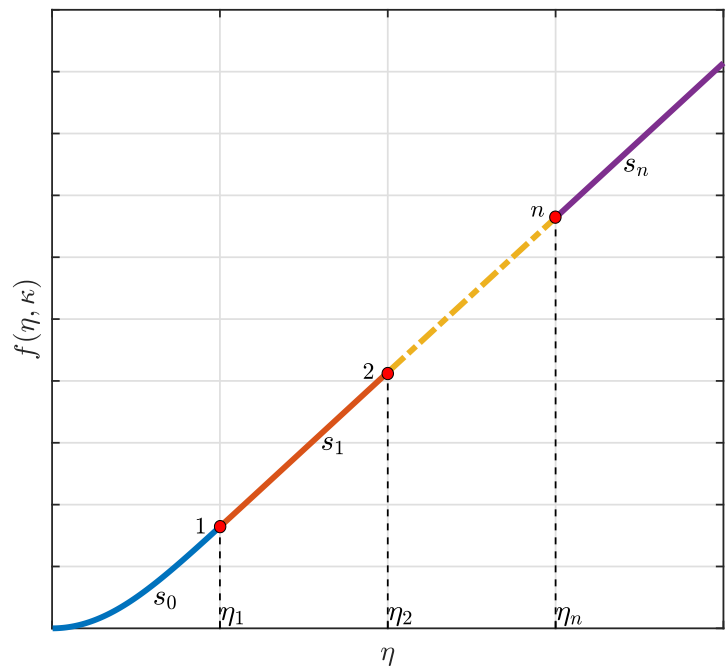}
    \caption{Taylor series expansions about leaping centers at  $\eta_{0} = 0, \eta_{1}, \eta_{2}, ...., \eta_{n}$. The symbol \(n\) denotes the number of expansions.}
    \label{Pic1}
\end{figure}

The idea used in the present approach is to generate a few Taylor series expansions for $f(\eta)$ with centers of expansion denoted as $\eta_{0} = 0, \eta_{1}, \eta_{2}, ...., \eta_{n}$, as illustrated in Figure(\ref{Pic1}). The corresponding Taylor's series expansions are denoted as $s_0,~s_1, \hdots, s_n$ for the functions \(f_{s_0},~f_{s_1}, \hdots, f_{s_n}\) , respectively. It may be noted that $f_{s_0}=f_{s_1}= \hdots=f_{s_n}$, because all the series expansions are for the same function, but the difference is only in their convergence properties.  We use Differential Transform method described in the above section to compute the coefficients of the Taylor's series. The coefficients of the series expansions are denoted by $F_{s_i}(k)$ for all the values of $k\in \mathbb{W}$. The initial condition for $s_{0}$ are  $F_{s_0}(0) = 0$, $F_{s_0}(1) = 0$ for no-slip flow and $F_{s_0}(1) = v_{slip}$, a prescribed values of slip flow. We keep $F_{s_0}(2) = \kappa^{*}$, a variable which is varied during the computation of Taylor's series coefficients. 

The initial conditions for computing any series with center of expansion at \(\eta_i\), namely,  $s_i$ are  taken from the already obtained series $s_{i-1}$ as,
\begin{equation}
    \begin{split}
        f_{s_i}(\eta_i) & = f_{s_{i-1}}(\eta_i)\\
        f_{s_i}'(\eta_i) & = f'_{s_{i-1}}(\eta_i)\\
        f_{s_i}''(\eta_i) & = f''_{s_{i-1}}(\eta_i)       
    \end{split}
\end{equation}
Note that $f_{s_{i-1}}(\eta)$ should converge at $\eta = \eta_{i}$. So fixing an $\eta_{i}$ for expanding a series $s_{i}$ depends on the convergence of $s_{i-1}$ at $\eta_i$. By following this procedure, it is possible to develop any number of series expansions. But an accurate  solution can be obtained from as low as three series expansions, which is going to be demonstrated in the section on results.

\subsection{Newton-Raphson Iteration}
We first develop the selected number of  function expansions,  $s_i$, \(i\in \mathbb{W}\) and  $i\le n$, by using a guessed value for \(\kappa^{*}\), say $\kappa^{*} = 0.1$. We then develop one more set of expansion namely, $s'_i$, corresponding to a slightly different value for \(\kappa^{*}\), denoted as $\kappa^{*}+\delta \kappa$, where $\delta \kappa$ is a very small value. A typical value for $\delta \kappa$ used in the present computations is $10^{-18}$. It is found that \(f'(\infty)=1\) is a strong boundary condition for the problem compared to \(f''(\infty)=0\). So we solve the equation  \(f'(\infty)-1=0\) for finding the converged solution for \(\kappa^{*}\). The converged solution of  \(\kappa^{*}\) is denoted as \(\kappa\). In the present computations the condition at infinity is changed to the one corresponding to a very high value of \(\eta\) denoted as \(\tilde{\infty}\). So we use \(f'(\tilde{\infty})-1=0\) as the equation solved for the determination of \(\kappa\). We have used a value for \(\tilde{\infty}\)  in the range \(5-30\), where $\tilde{\infty}$ is a finite value of $\eta$ at which $f'(\eta)$ becomes almost equal to 1. The approximate gradient for Newton-Raphson method is computed as,
\begin{equation}
    g = \frac{f'_{s'_{n}}(\tilde{\infty}) - f'_{s_{n}}(\tilde{\infty})}{\delta \kappa}
\end{equation}
$\kappa^{*}$ is now updated as,
\begin{equation}
    \kappa^{*}_{new} = \kappa^{*} - \frac{(f'(\tilde{\infty}) - 1.)}{g}
\end{equation}
Now, $\kappa^{*}_{new}$ is assigned as the new value of \(\kappa^{*}\) and the computation of series is continued  till $f'(\tilde{\infty})$ becomes equal to 1 and this stage the value of \(\kappa^{*}\) is treated as the converged  value for \(\kappa\). The number of Newton-Raphson iterations applied to reach the solution for \(\kappa\) is kept as \(50\).
\section{Results and discussion} 
The methodology for the solution has been implemented in FORTRAN with variable type for computations set using "IMPLICIT REAL*16(A-H,O-Z)" so as to get accuracy up to 32 decimal places in the floating point operations. All the results we report here are obtained from summing \(2000\) terms of the Taylor's series.  We find no difference in the converged value of \(\kappa\) with further increase in the number of terms of Taylor's series summed.
\subsection{Flow with no-slip}
\begin{figure}[htp!]
    \centering
    \includegraphics[scale=0.55]{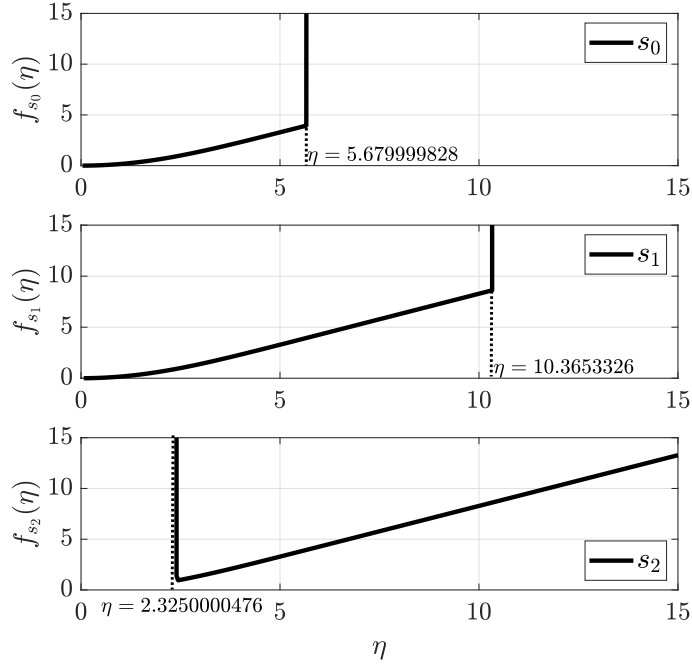}
    \caption{Convergence properties of  three series expansions, $f_{s_{0}}$, $f_{s_{1}}$ and $f_{s_{2}}$ having centers of expansion  \(\eta_0=0\), \(\eta_1=5\) and \(\eta_2=10\).}
    \label{conv-pro}
\end{figure}

\begin{figure}[htp!]
     \centering
     \begin{subfigure}[b]{0.485\textwidth}
         \centering
         \includegraphics[width=\textwidth]{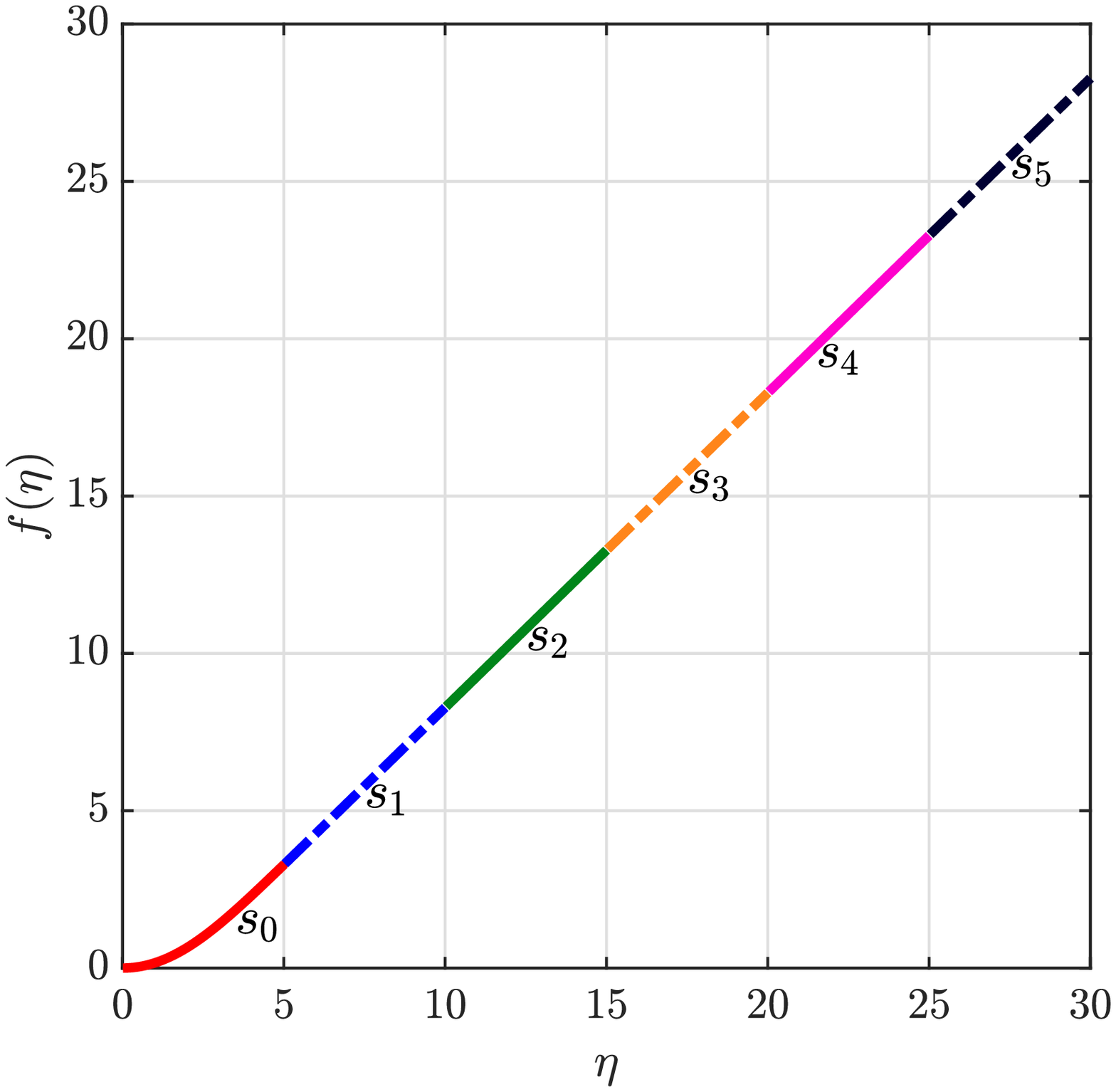}
         \caption{$f(\eta)$ vs. $\eta$}
         \label{f}
     \end{subfigure}
     \hfill
     \begin{subfigure}[b]{0.485\textwidth}
         \centering
         \includegraphics[width=\textwidth]{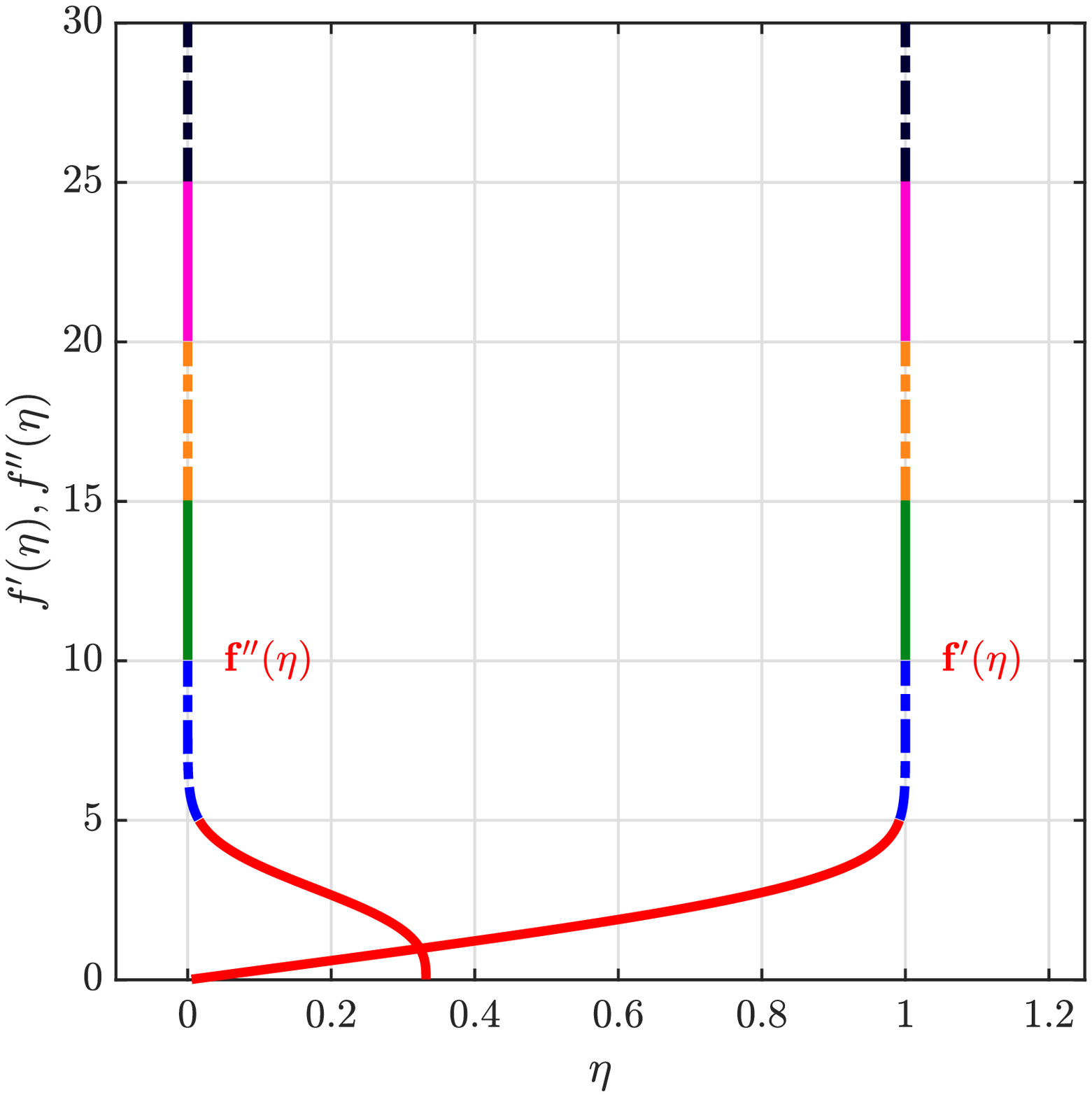}
         \caption{$f^{'}(\eta)$ and $f^{''}(\eta)$ vs. $\eta$}
         \label{fdAndfdd}
     \end{subfigure}
        \caption{The accurate solutions (a)  $f(\eta)$, (b) $f'(\eta)$ and $f''(\eta)$  for $\eta\in[0,30]$.}
        \label{SuccessOfNewMethod}
\end{figure}

\begin{table}[htp!]
    \centering
    \begin{tabular}{c|c|c} 
    \hline
    n  & \(\tilde{\infty}\) & \(f''(0)\)  \\
    \hline
     1  & 5  &  0.3361523439045316853836761393072\\
       2& 10  &  0.3320573372037024527094809622944 \\
        3 & 15 &  0.33205733621519629893842022877817 \\
         4& 20 & 0.33205733621519629893718006201058 \\
         5& 25 & 0.33205733621519629893718006201058 \\
6 & 30 & 0.33205733621519629893718006201058 \\
\hline
    \end{tabular}
    \caption{Effect of number of expansion $n$ on the converged solution of the value of \(f''(0)\).}
    \label{effect}
\end{table}
Table \ref{effect} shows the effect of number of expansions \(n\) on the converged solution of \(\kappa\) by leaping the center of expansion by  \(\triangle\eta=5\). Which means for every \(s_i\), the center of expansion \(\eta_i=\eta_{i-1}+5\), with \(\eta_0=0\).  The number of iterations of Newton-Raphson step  is kept as \(50\) and this is found to guarantee the convergence of \(f''(0)\), up to \(32\) decimal places. It can be seen that the converged solution is independent of \(n\) for \(n\ge 4\). This implies that the boundary condition at infinity can be well approximated at \(\tilde{\infty}=20\). More precisely, we find that the value of \(\kappa\) does not change from \(\tilde{\infty}=18\) onward.   We have also examined  the effect of change of separation between the center of expansions and the observations are tabulated in Table \ref{expansionEffect}. It can be seen that the result is independent on the step size of leaping of the center of expansion and is only dependent on the point \( \tilde{\infty}\), where the boundary condition \(f'(\eta)=1\) is applied. 

 The value of \(\kappa\) resulted from the present computation has 32 decimal points accuracy, the accuracy is better than that of 0.33205733621519630 as given by Boyd in \cite{Boyd2008} and 0.332057336215196299183 as given by Lal and Neeraj in \cite{Lal2014}.
\begin{table}[htp!]
    \centering
    \begin{tabular}{c|c|c|c} 
    \hline
   \(\eta_i-\eta_{i-1}\) &  n  & \(\tilde{\infty}\) & \(f''(0)\)  \\
    \hline
2  &    10  & 20  & 0.33205733621519629893718006201058\\
3 &      7 & 21  &  0.33205733621519629893718006201058 \\
4 &      5 & 20 & 0.33205733621519629893718006201058 \\
5 &     4    & 20 & 0.33205733621519629893718006201058 \\
\hline
    \end{tabular}
    \caption{Effect of change of leaping distance of consecutive centers of expansion}
    \label{expansionEffect}
\end{table}

The convergence properties of the series expansions for \(s_0\), \(s_1\) and \(s_2\) having centers of expansion at \(\eta_0=0\), \(\eta_0=5\) and \(\eta_0=10\) are depicted in Figure \ref{conv-pro} for \(\eta\in[0,15]\). It can be seen that the series \(s_0\) converges smoothly up to \(\eta=5.69\), which is close to the value corresponding to the singularity at \(\eta=-5.69\) as reported by Boyd in \cite{Boyd2008}. The series \(s_1\) is found to be converging for \(\eta\in[0,10.365]\) and diverges for \(\eta>10.365\). The series \(s_2\) is converging for \(\eta\in[2.4,15]\) and diverges for \(\eta<2.4\). The radius of convergence increases with the increase in distance of center of expansion from \(\eta=0\). These observations on the convergence implies that a single series can not represent the  expansion of Blasius function for a larger domain.    

By considering the convergence properties, we propose to use a combination of leaping series expansions to represent Blasius function as demonstrated in  Figure (\ref{SuccessOfNewMethod}). It consists of \(6\) expansions, namely \(s_0\) to \(s_5\) with centers of expansion at \(\eta_0=0\), \(\eta_1=5,~\hdots\), \(\eta_i=\eta_{i-1}+5\) for determining the convergence properties upto $\eta = 30$. The series expansion \(s_i\) is used to determine the value of the function in the converging interval,  \(\eta\in[\eta_i,\eta_i+5]\). Figure (\ref{SuccessOfNewMethod}) demonstrates the expansion of the function in the interval \(\eta\in[0,30]\). In the same manner, by combining more number of series expansions, it is possible to expand the function for any positive real value. As per the results obtained in the present study \(\eta=20\) found to be an optimum point  where the boundary condition \(f'(\infty)=1\)  can be applied. From \(\eta=18\) onward the asymptotic variation of the function is as \(f(\eta)= -1.720787657520502819605438159825+\eta\). The parameter in the asymptotic variation has more decimal points accuracy than the one reported by Boyd in \cite{Boyd2008}. Tables \ref{tableFandFd} and \ref{tablefddFddd} gives the benchmark results of the present accurate solutions of \(f(\eta)\), \(f'(\eta)\), \(f''(\eta)\) and \(f('''\eta)\) for \(\eta\in[0,15]\). It can be seen that the data  satisfy the equation \(f'''+ff''/2=0\) almost accurately.  
 
\begin{table}[htp!]
\caption{Benchmark solutions of \(f(\eta)\) and \(f'(\eta)\) for \(\eta\in[0,15]\)  of no-slip Blasius boundary layer.}
\label{tableFandFd}
\centering
\resizebox{\textwidth}{!}{\begin{tabular}{c|c|c}
\hline
\(\eta\) & \(f(\eta)\)                         & \(f'(\eta)\)                       \\ \hline
0.       & 0.00000000000000000000000000000000  & 0.00000000000000000000000000000000 \\
1.       & 0.16557172578927971994607139866464  & 0.32978003124966696806286485058647 \\
2.       & 0.65002436993528859325840008060347  & 0.62976573650238585970789214001077 \\
3.       & 1.39680823087034546785665318114741  & 0.84604444365799349725098240374192 \\
4.       & 2.30574641846207117902246203651341  & 0.95551822981069424686181417757482 \\
5.       & 3.28327366515631417796776890035073  & 0.99154190016439334277083771840218 \\
6.       & 4.27962092251384908051614632858349  & 0.99897287243586052111975290945127 \\
7.       & 5.27923881102910880195475622555865  & 0.99992160414795010955867185310631 \\
8.       & 6.27921343134607434184668676342406  & 0.99999627453530080913338909825509 \\
9.       & 7.27921237064524251498448401453620  & 0.99999989045379928945150717872348 \\
10.      & 8.27921234293432582800203409779852  & 0.99999999801539054606962551624450 \\
11.      & 9.27921234248405790426225183944242  & 0.99999999997791640300822743404095 \\
12.      & 10.27921234247952546306776851892884 & 0.99999999999984941269916627913591 \\
13.      & 11.27921234247949728852564004547502 & 0.99999999999999937182753490370178 \\
14.      & 12.27921234247949718064881057910890 & 0.99999999999999999839911378268159 \\
15.      & 13.27921234247949718039492878137463 & 0.99999999999999999999751013528420 \\ \hline
\end{tabular}}
\end{table}

\begin{table}[htp!]
\caption{Benchmark solutions of \(f''(\eta)\) and \(f'''(\eta)\) for \(\eta\in[0,15]\) of no-slip Blasius boundary layer.}
\label{tablefddFddd}
\centering
\resizebox{\textwidth}{!}{\begin{tabular}{c|c|c}
\hline
\(\eta\) & \(f''(\eta)\)                      & \(f'''(\eta)\)                      \\ \hline
0.       & 0.33205733621519629893718006201058 & 0.00000000000000000000000000000000  \\
1.       & 0.32300711668694284746919034374890 & -0.02674042287603818941557199292502 \\
2.       & 0.26675154569727845645745365043572 & -0.08669750271056888580038083842490 \\
3.       & 0.16136031954087853155687842360081 & -0.11269471123528408855663373868942 \\
4.       & 0.06423412109169066819653982932531 & -0.07405379732511237191225366739570 \\
5.       & 0.01590679868531816315507141064886 & -0.02611318661022410269690678475316 \\
6.       & 0.00240203984375727885991930092577 & -0.00513990998602777383188838650274 \\
7.       & 0.00022016895527113462513633196741 & -0.00058116224682555289777061646404 \\
8.       & 0.00001224092624325319054622340463 & -0.00003843169423937603885233795627 \\
9.       & 0.00000041278790156878232013889098 & -0.00000150238539977608548717786204 \\
10.      & 0.00000000844291586701753292682574 & -0.00000003495034662828381184818483 \\
11.      & 0.00000000010473955482006946773795 & -0.00000000048595028491633710186547 \\
12.      & 0.00000000000078810075727866872385 & -0.00000000000405052751566817613004 \\
13.      & 0.00000000000000359671084986261923 & -0.00000000000002028403270505018845 \\
14.      & 0.00000000000000000995593526876784 & -0.00000000000000006112552161659052 \\
15.      & 0.00000000000000000001671519173257 & -0.00000000000000000011098229018106 \\ \hline
\end{tabular}}
\end{table}

\subsection{Flow with slip}

\begin{figure}[htp!]
     \centering
     \begin{subfigure}[b]{0.485\textwidth}
         \centering
         \includegraphics[width=\textwidth]{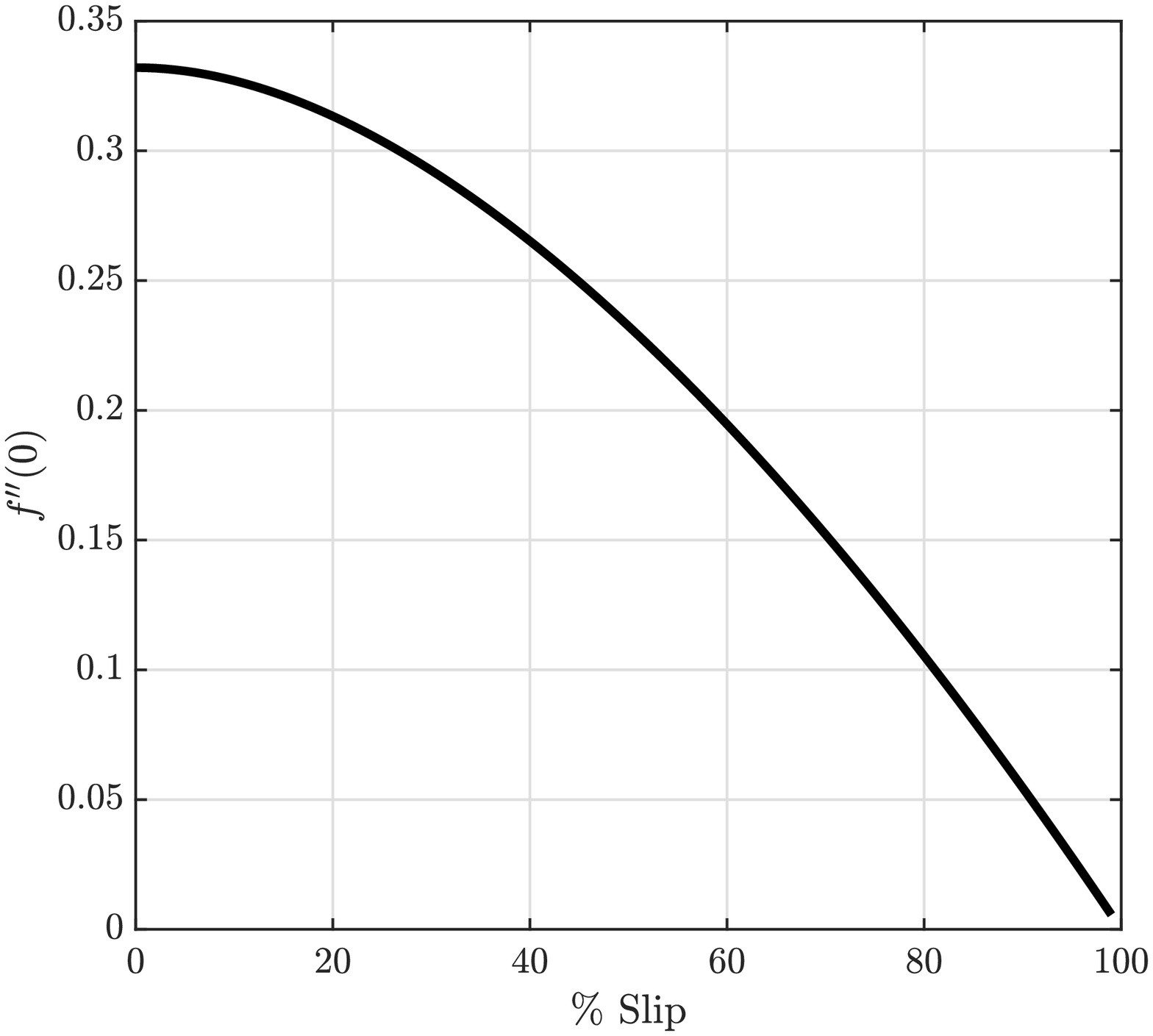}
         \caption{\(f''(0)\) Vs \(\%\) slip}
         \label{SlipFdd}
     \end{subfigure}
     \hfill
     \begin{subfigure}[b]{0.485\textwidth}
         \centering
         \includegraphics[width=\textwidth]{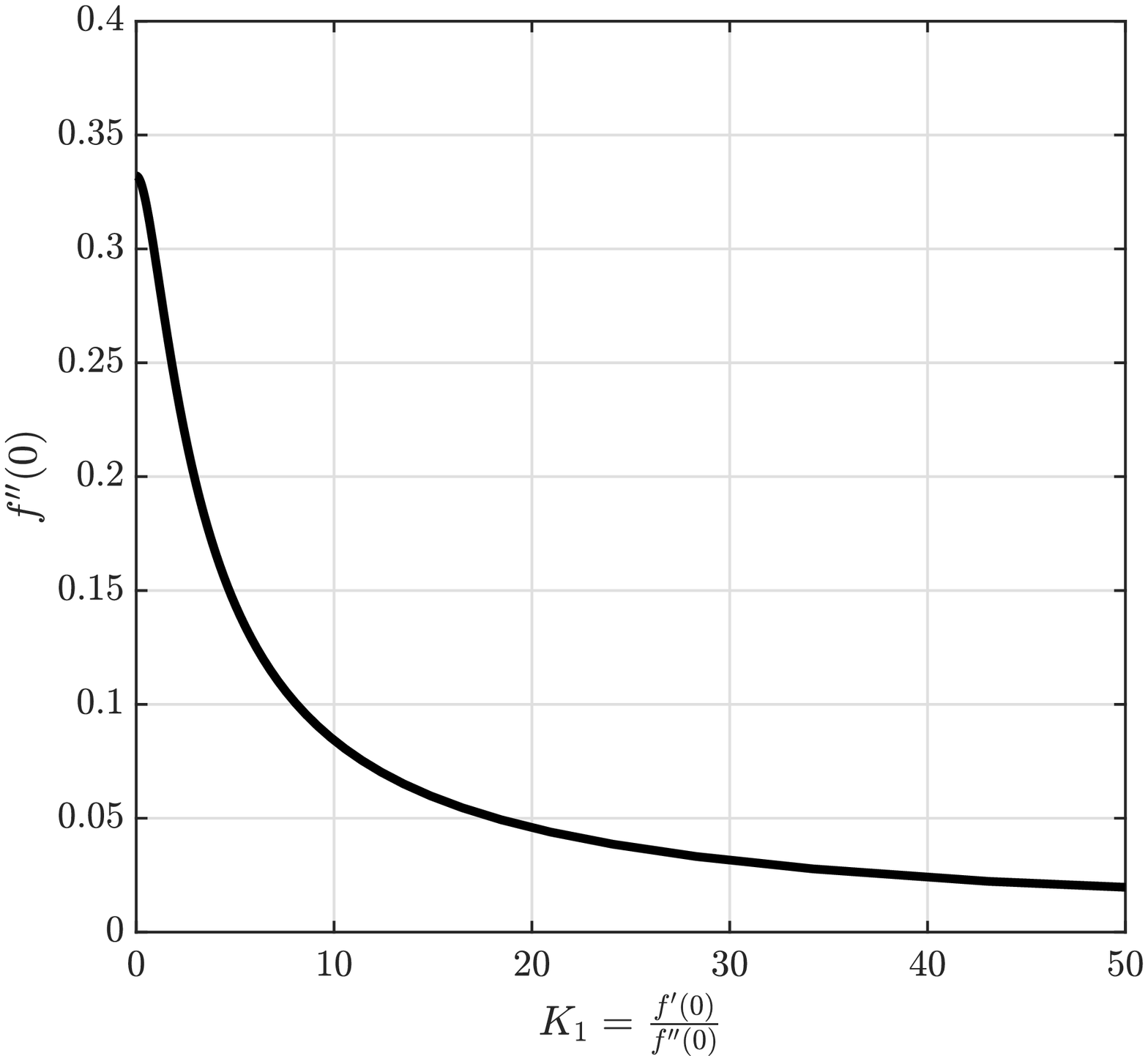}
         \caption{ $f''(0)$ Vs \(K_1\)}
         \label{SlipK1}
     \end{subfigure}
        \caption{Variation of \(f''(0)\) with \(\%\)slip and the rarefaction parameter \(K_1\).}
        \label{SlipGraphs}
\end{figure}

\begin{figure}[htp!]
     \centering
     \begin{subfigure}[b]{0.485\textwidth}
         \centering
         \includegraphics[width=\textwidth]{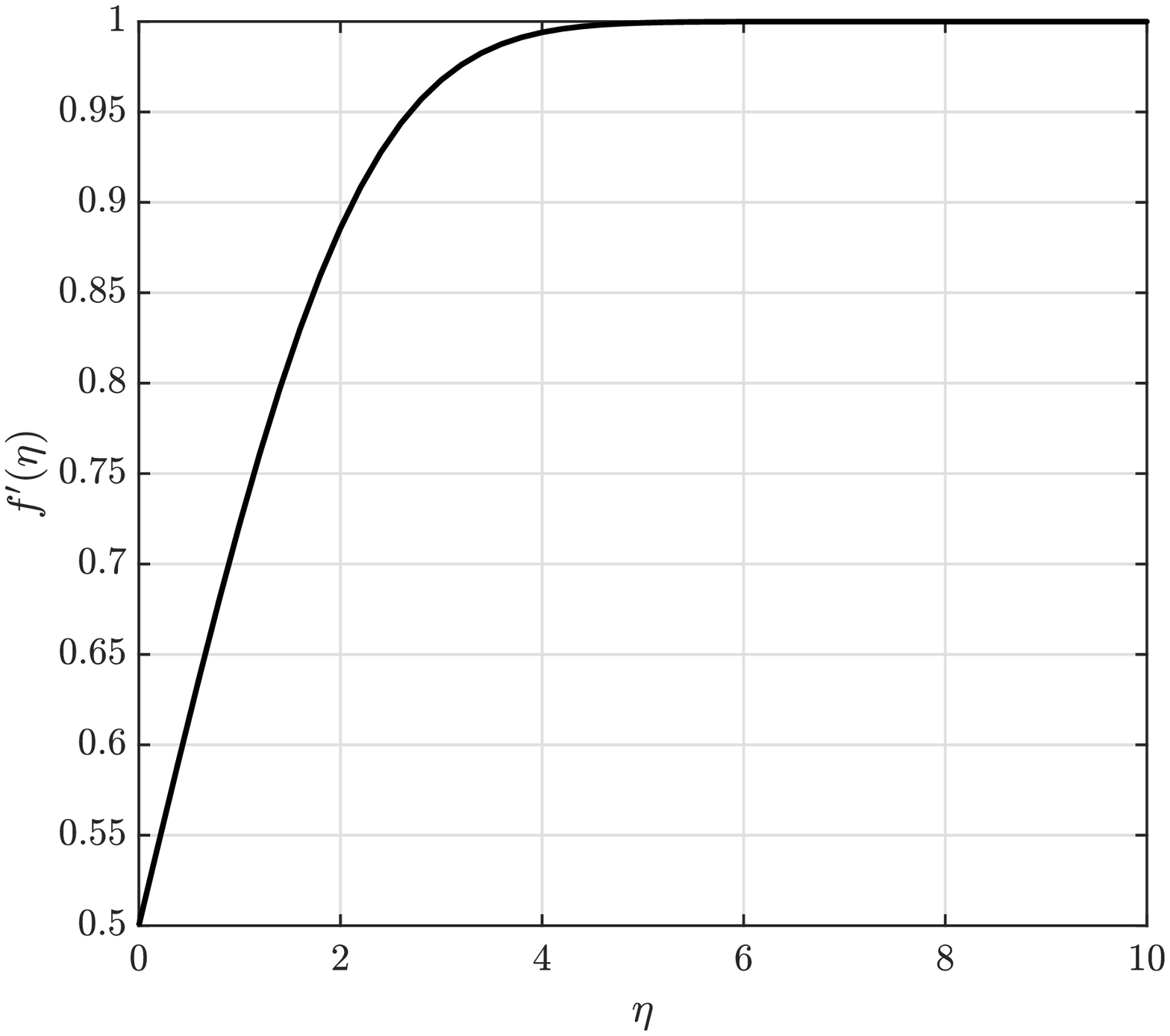}
         \caption{\(f'(\eta)\) Vs $\eta$}
         \label{SlipF1d5}
     \end{subfigure}
     \hfill
     \begin{subfigure}[b]{0.485\textwidth}
         \centering
         \includegraphics[width=\textwidth]{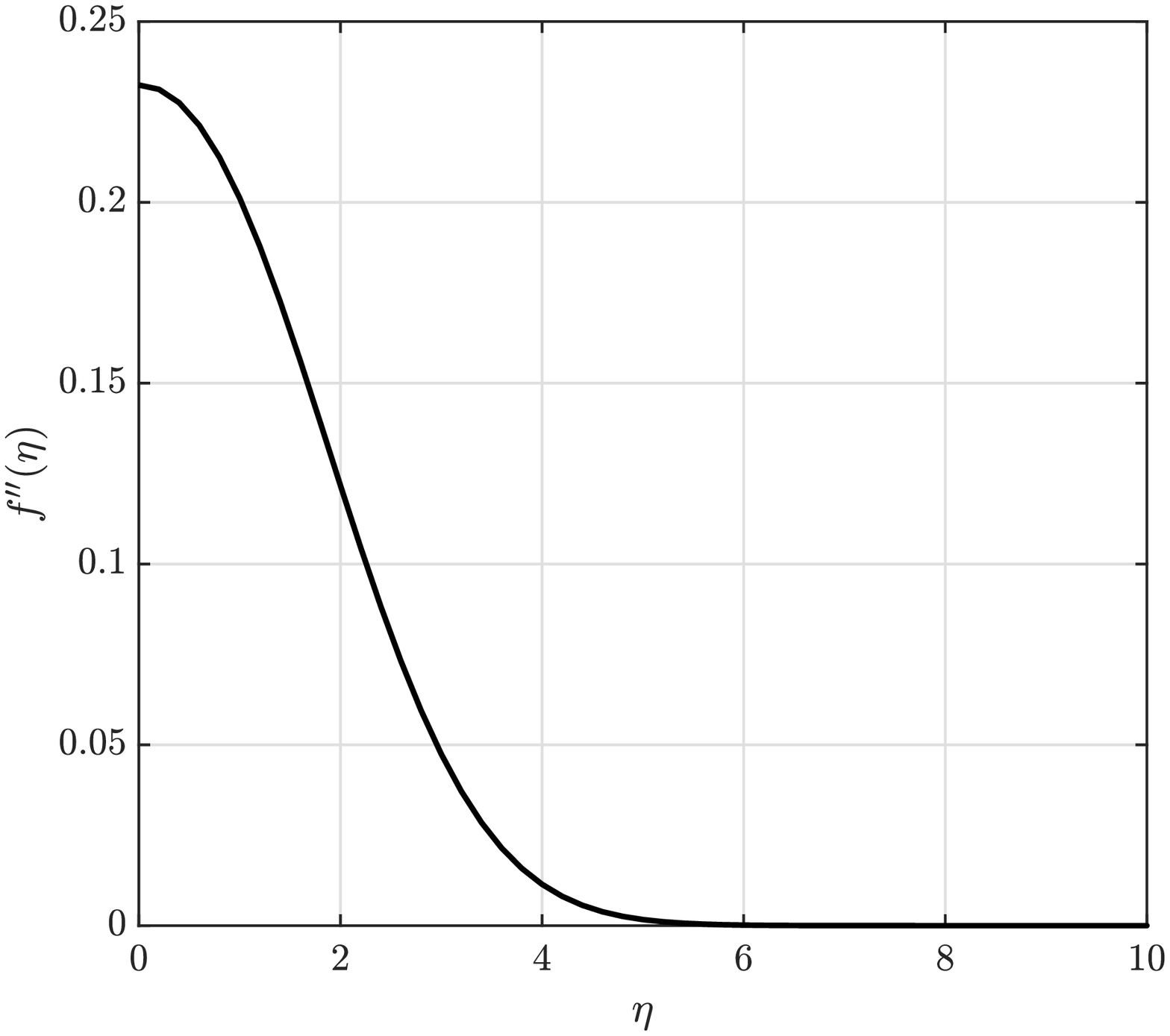}
         \caption{ \(f''(\eta)\) Vs $\eta$}
         \label{SlipF1dd5}
     \end{subfigure}
        \caption{$f'(\eta)$ and $f''(\eta)$ plotted for different values of $\eta$ when the slip is 0.5}
        \label{slip5}
\end{figure}

The present methodology for computing \(f''(0)\) using Newton-Raphson method has enabled computation of the same for the cases with velocity slip  boundary condition at the wall. We have varied the slip   velocity, $v_{slip}\in [0,1]$, with the extreme values \(0\) and \(1\) imply  no slip and uniform velocity, respectively. Figure \ref{SlipFdd} shows the plot showing the variation of  \(f''(0)\) with the percentage slip. The curve has a monotonic decreasing non-linear variation. It can be found that  \cite{Martin2001} has considered the ratio of slip to wall shear stress as $K_{1}$,  which represents the rarefaction of the flow medium, surface properties etc., as given in equation(\ref{surf}). We have computed the solution corresponding to \(K_1\in[0,50]\) and the result is presented in Figure \ref{SlipK1} as the variation of $f''(0)$ with $K_{1}$. It is evident that with increase of $K_{1}$, the value of \(f''(0)\) reduces and  slope of the variation reduces continuously and is expected to reach zero asymptotically. 
\begin{equation}
    K_{1} = \frac{2-\sigma}{\sigma}Kn_{x}{Re_{x}}^{1/2} \label{surf}
\end{equation}
Where, $\sigma$ is called  tangential accommodation coefficient,  \(Kn_x\) and \(Re_x\) are the local Knudsen and Reyonld numbers. Table (\ref{benchmarkSlip}) provides benchmark numerical values  of \(f''(0)\) for five typical values of the slip. Figure \ref{slip5} shows the variations of \(f'(\eta)\) and \(f''(\eta)\) for the case with slip\(=0.5\).

\section{Conclusions}
In order to overcome the convergence problems of a single series expansion of Blasius function, a combination of series expansions with different convergence properties, that are expanded about different centers called the leaping centers is demonstrated.  Though development of series expansion of functions governed by differential equations requires all the required conditions as initial values, we demonstrated deduction of series solutions for boundary value problems using Newton-Raphson method in this work. The major conclusions from this study are:
\begin{enumerate}
    \item The solution is found to be independent on the point of application of the boundary condition \(f'(\infty=1.)\) for \(\eta\ge 18\). 
    \item The series with the center of expansion at \(\eta=0\) diverges for \(\eta>5.68\). The series with the center of expansion at \(\eta=5\) diverges for \(\eta>10.36\). The series with the center of expansion at \(\eta=10\) has found to exhibit divergence for \(\eta<2.32\).
    \item The solution is found to be independent of the leaping distance between the center of expansions for leaping distances less than \(5.68\), corresponding to the divergence of the series with the center of expansion at \(\eta=0\).
    \item The variation of \(f''(0)\) with boundary slip is found to be monotonically decreasing. The variation of \(f''(0)\) with the rarefaction parameter has a decreasing tend having an asymptotic behaviour.
    \item The accuracy of the variations and all the related parameters of Blasius function are obtained with a better accuracy than those atready available in the literature.
    \item The method can be used for solving differential equations of functions having a singularity on the left side of the origin.
    
\end{enumerate}

\begin{table}[htp!]
\caption{Variation of $f''(0)$ with wall slip condition.}
\label{benchmarkSlip}
\centering
\begin{tabular}{c|c}
\hline
Slip & \(f''(0)\)                             \\ \hline
0.2      & 0.31335775171380422946976167600464 \\
0.4      & 0.26523233119785271052478203825874 \\
 0.5 &   0.23245507775976810752799593194116 \\
0.6      & 0.19462617044304093400347255410066 \\ 
0.8      & 0.10538860552764686718424828804120 \\\hline
\end{tabular}
\end{table}

\bibliographystyle{elsarticle-num-names} 
\bibliography{ref}





\end{document}